\documentclass[11pt]{article}
\usepackage{amsmath,amsfonts,amsthm,enumerate}
\usepackage[colorlinks=true,citecolor=black,linkcolor=black,urlcolor=blue]{hyperref}
\usepackage[lmargin=38mm,rmargin=38mm,tmargin=30mm,bmargin=30mm]{geometry}
\usepackage[sort&compress,numbers]{natbib}

\renewcommand{\baselinestretch}{1.1}
\renewcommand{\thefootnote}{\fnsymbol{footnote}}	
\allowdisplaybreaks

\newcommand{\arXiv}[1]{arXiv:\,\href{http://arxiv.org/abs/#1}{#1}}
\newcommand{\msn}[1]{MR:\,\href{http://www.ams.org/mathscinet-getitem?mr=MR#1}{#1}}
\newcommand{\MSN}[2]{MR:\,\href{http://www.ams.org/mathscinet-getitem?mr=MR#1}{#1}}
\newcommand{\Zbl}[1]{Zbl:\,\href{http://www.zentralblatt-math.org/zmath/en/search/?q=an:#1}{#1}}
\newcommand{\doi}[1]{doi:\,\href{http://dx.doi.org/#1}{#1}}

\newcommand{\urlprefix}{}

\theoremstyle{plain}
\newtheorem{theorem}{Theorem}
\newtheorem{lemma}[theorem]{Lemma}
\newtheorem{corollary}[theorem]{Corollary}

\theoremstyle{definition}

\newcommand{\thmlabel}[1]{\label{thm:#1}}
\newcommand{\thmref}[1]{Theorem~\ref{thm:#1}}
\newcommand{\lemlabel}[1]{\label{lem:#1}}
\newcommand{\lemref}[1]{Lemma~\ref{lem:#1}}
\newcommand{\twolemref}[2]{Lemmas~\ref{lem:#1} and \ref{lem:#2}}
\newcommand{\corlabel}[1]{\label{cor:#1}}
\newcommand{\corref}[1]{Corollary~\ref{cor:#1}}

\def\bag{\mathcal{B}}

\begin{document}

\title{\bf Cliques in Odd-Minor-Free Graphs\,\footnotemark[1]}
\author{Ken-ichi Kawarabayashi\,\footnotemark[2] \qquad David  R.~Wood\,\footnotemark[3]}
\footnotetext[1]{A preliminary version of this paper was published in the Proceedings of Computing: the Australasian Theory   Symposium (CATS 2012).}
\footnotetext[2]{National Institute of Informatics, Tokyo, Japan (\texttt{k\_keniti@nii.ac.jp}).}
\footnotetext[3]{Department of Mathematics and Statistics, The University of
Melbourne, Melbourne, Australia
(\texttt{woodd@unimelb.edu.au}). Supported by a QEII Fellowship and
Discovery Project from the Australian Research Council.}

\maketitle

\begin{abstract}
This paper is about: (1) bounds on the number of 
cliques in a graph in a particular class, and (2) algorithms for
listing all cliques in a graph. We present a simple algorithm that
lists all cliques in an $n$-vertex graph in $O(n)$ time per clique. 
For $O(1)$-degenerate graphs, such as graphs excluding a fixed minor,
we describe a $O(n)$ time algorithm for listing all cliques. 
We prove that graphs excluding a fixed odd-minor have $O(n^2)$
cliques (which is tight), and conclude a $O(n^3)$ time algorithm for
listing all cliques. 
\end{abstract}

\section{Introduction}

\renewcommand{\thefootnote}{\arabic{footnote}}	

A \emph{clique} in a graph\footnote{We consider simple finite
  undirected graphs $G$ with vertex set $V(G)$ and edge set
  $E(G)$. For each vertex $v\in V(G)$, let $N_G(v)$ or simply $N(v)$,
  be $\{w\in V(G):vw\in E(G)\}$.}
is a set of pairwise adjacent vertices.  This paper is about:
\begin{enumerate}[(1)]
\item bounds on the number of
cliques in a graph in a particular class of graphs, and
\item algorithms for listing all cliques in a graph in
such a class.
\end{enumerate}
In addition to being of intrinsic interest, bounds on the number of
cliques in a graph have recently been used in a proof that
minor-closed graph classes are `small' \citep{NSTW-JCTB06}, and in the
analysis of a linear-time algorithm for computing separators in graphs
in minor-closed classes \citep{ReedWood-TALG}, which in turn has been
applied in shortest path \citep{WN10,MT-DAM09,Yuster-TCS} and maximum
matching \citep{YZ-SODA07} algorithms. Note that (1) and (2) for
\emph{maximal} cliques have been extensively studied; see \citep{ELS-ISAAC}
and the references therein.

This paper describes 
a simple algorithm that lists all cliques in a given $n$-vertex graph in $O(n)$ time per clique
(\thmref{NonRecCliques}). This implies that if we solve (1) for a
particular class, then we immediately solve (2). Note that analogous
results hold for maximal cliques: there are algorithms that list all
maximal cliques in polynomial time per clique
\citep{Eppstein-TAlg09,GNS-DAM09,LLR-SJC80,TIS-SJC77,JYP-IPL88,PU59} or in total time
proportional to the maximum possible number of cliques in an
$n$-vertex graph, without additional polynomial factors
\citep{Eppstein-JGAA03,TTT-TCS06}.

As an example of (1), many authors have observed that every $n$-vertex
planar graph contains $O(n)$ cliques
\citep{Eppstein-JGT93,PY-IPL81}. \citet{Wood-GC07} proved the best
possible upper bound of $8(n-2)$. More generally, for each surface
$\Sigma$, \citet{DFJW} characterised the $n$-vertex graphs embeddable
in $\Sigma$ with the maximum number of cliques in terms of so-called
irreducible triangulations. They also proved that if $K_{\omega}$ is
the largest complete graph that embeds in $\Sigma$, then every
$n$-vertex graph that embeds in $\Sigma$ contains at most
$8n+\frac{3}{2}2^{\omega}+o(2^{\omega})$ cliques. Exact results and a
precise characterisation of the extremal examples are obtained for
graphs that embed in the plane, torus, double torus, projective plane,
$\mathbb{N}_3$, and $\mathbb{N}_4$.

These results are generalised by considering $H$-minor-free graphs.  A
graph $H$ is a \emph{minor} of a graph $G$ if $H$ can be obtained from
a subgraph of $G$ by contracting edges. Equivalently, $H$ is a minor
of $G$ if $G$ contains a set of vertex-disjoint trees, one tree $T_v$
for each vertex $v$ of $H$, such that for every edge $e=vw$ in $H$
there is an edge $\hat e$ between $T_v$ and $T_w$. A graph $H$ is an
\emph{odd minor} of $G$ if, in addition, the vertices in
$\bigcup_vV(T_v)$ can be 2-coloured such that for each vertex $v\in
V(H)$ the edges in $T_v$ are bichromatic, and for each edge $e=vw\in
E(H)$, the edge $\hat e$ between $T_v$ and $T_w$ is monochromatic. A
graph is (\emph{odd-})\emph{$H$-minor-free} if it contains no
(odd-)$H$-minor.

Several authors have proved that for every fixed graph $H$, every
$H$-minor-free graph with $n$ vertices contains $O(n)$ cliques
\citep{FOT,Wood-GC07,ReedWood-TALG,NSTW-JCTB06}. The best bound, due
to \citet{FOT}, states that every $K_t$-minor-free graph contains at
most $c^{t\log\log t}n$ cliques, for some constant $c$. It is open
whether such graphs have at most $c^tn$ cliques \citep{Wood-GC07}.

This paper considers (1) and (2) for graphs that exclude an odd
minor. The class of odd-$H$-minor-free graphs is more general than the
class of $H$-minor-free graphs.  For example, the complete bipartite
graph $K_{n,n}$ contains a $K_{n+1}$ minor but contains no
odd-$K_3$-minor. In fact, a graph contains no odd $K_3$-minor if and
only if it is bipartite. In general, every $K_t$-minor-free graph has
$O(t\sqrt{\log t} n)$ edges, and this bound is best possible
\citep{Kostochka84,Thomason84}. On the other hand, some
odd-$K_t$-minor-free graphs, such as $K_{n,n}$, have $\Theta(n^2)$
edges. This paper proves the following theorem:

\begin{theorem}
  \thmlabel{OddCliques} For every fixed graph $H$, there is a constant
  $c$, such that every $n$-vertex odd-$H$-minor-free graph $G$
  contains at most $cn^2$ cliques, and these cliques can be listed in
  $O(n^3)$ time.
\end{theorem}

The bound on the number of cliques in \thmref{OddCliques} is best
possible up to the value of $c$, since $K_{n,n}$ contains no
odd-$K_3$-minor and contains $\Theta(n^2)$ cliques. Also note that a
polynomial bound on the number of cliques in every graph in a class is
non-trivial, since $K_n$ contains $2^n$ cliques.

\thmref{OddCliques} is in sharp contrast with a number of
intractability results about finding cliques: it is NP-complete to
test if a graph $G$ contains a $k$-clique (given $G$ and $k$)
\citep{Karp72}; it is $W[1]$-complete to test if a graph $G$ contains
a $k$-clique (given $G$ with parameter $k$) \citep{DF-TCS95}; and
approximating the maximum clique size is hard \citep{ALMSS}.


\section{General Graphs}

Consider the following simple recursive algorithm for listing all
cliques in a graph.
\begin{center}
\framebox{
  \begin{minipage}{120mm}
    \textsc{Cliques}($G$)\\
    \emph{input}: graph $G$\\
    \emph{output}: the set of all cliques in $G$ \medskip \hrule
    \medskip
    \begin{tabular}{ll}
      1. & if $V(G)=\emptyset$ then return $\{\emptyset\}$\\
      2. &  choose $v\in V(G)$\\
      3. &   return $\{C\cup\{v\}: C\in  \textsc{Cliques}(G[N_G(v)])\,\}\bigcup\textsc{Cliques}(G-v)$
    \end{tabular}
  \end{minipage}} 
\end{center}
\medskip

\begin{theorem}
  \thmlabel{GenCliques} If $G$ is an $n$-vertex graph then
  $\textsc{Cliques}(G)$ returns the set of all cliques in $G$.
\end{theorem}

\begin{proof}
  We proceed by induction on $|V(G)|$. If $V(G)=\emptyset$ then
  $\emptyset$ is the only clique in $G$, and the algorithm correctly
  returns the set of all cliques in $G$.  Otherwise, each clique $C$
  of $G$ either contains $v$ or does not contain $v$. In the first
  case, $C$ is a clique of $G$ containing $v$ if and only if
  $C=S\cup\{v\}$ for some clique $S$ of $G[N_G(v)]$.  In the second
  case, $C$ is a clique of $G$ not containing $v$ if and only if $C$
  is a clique of $G-v$.  Therefore, by induction, the algorithm
  correctly returns the set of all cliques of $G$.
\end{proof}

The next algorithm outputs all cliques in $O(n)$ time per clique.
\begin{center}
\framebox{
  \begin{minipage}{120mm}
    \textsc{AllCliques}($G$)\\
    \emph{input}: graph $G$\\
    \emph{output}: all cliques in $G$ \medskip \hrule \medskip
    \begin{tabular}{ll}
      1. & output $\emptyset$\\
      2. & $i:=1$\\
      3. & $V_i:=V(G)$\\
      4. & repeat \\
      5. &    \hspace*{10mm} if $V_i=\emptyset$ then $i:=i-1$\\
      6. &   \hspace*{10mm} else\\
      7. &    \hspace*{20mm} choose $x_i\in V_i$\\
      8. &    \hspace*{20mm} output $\{x_1,\dots,x_i\}$\\
      9. &    \hspace*{20mm} $V_{i+1}:=V_i\cap N_G(x_i)$\\
      10. &    \hspace*{20mm} $V_i:=V_i\setminus\{x_i\}$\\
      11. &    \hspace*{20mm} $i:=i+1$\\
      12. &     \hspace*{10mm} end-if\\
      13. & until $i=0$
    \end{tabular}
  \end{minipage}}
\end{center}
\medskip

\begin{theorem}
  \thmlabel{NonRecCliques} If $G$ is a graph with $n$ vertices, 
then $\textsc{AllCliques}(G)$ outputs all cliques in $G$ in $O(n)$ time per clique.
\end{theorem}

\begin{proof}
  It is easily seen that \textsc{AllCliques} is simply a non-recursive
  implementation of \textsc{Cliques}, and therefore correctly
  outputs all cliques in $G$. To implement this algorithm efficiently,
  without loss of generality, assume that $V(G)=\{1,2,\dots,n\}$, and
  the adjacency lists and the sets $V_i$ are sorted.  Thus lines 7--11
  can be implemented in $O(n)$ time, and line 5 can be computed in
  $O(1)$ time.  Between outputting successive cliques, lines 7--11 are
  executed once, and line 5 is executed at most $n$ times. Thus the
  algorithm takes $O(n)$ time between outputting successive cliques.
\end{proof}

\section{Degenerate Graphs}

A graph $G$ is \emph{$d$-degenerate} if every non-empty subgraph of
$G$ has a vertex of degree at most $d$. For example, every planar
graph is 5-degenerate, and every $K_t$-minor-free graph is
$O(t\sqrt{\log t})$-degenerate \citep{Kostochka84,Thomason84}.
\citet{Wood-GC07} proved that every $d$-degenerate graph contains at
most $2^d(n-d+1)$ cliques, and this bound is tight for a $d$-tree.
Below we give an algorithm for finding all cliques in a $d$-degenerate
graph.

First consider the following data structure.  A linear ordering
$(v_1,\dots,v_n)$ of the vertices of a graph $G$ is
\emph{$d$-degenerate} if $|N^+(v_i)|\leq d$ for each vertex $v_i$,
where $N^+(v_i):=\{v_j:i<j,\,v_iv_j\in E(G)\}$. It is easily seen that
a graph is $d$-degenerate if and only if it has a $d$-degenerate
vertex ordering \citep{LW-CJM70}. Moreover, there are $O(dn)$ time
algorithms for computing a $d$-degenerate ordering of a given
$d$-degenerate graph, along with the set $N^+(v_i)$; see
\citep{CE91,ReedWood-TALG}. Also note that given a $d$-degenerate
ordering and given the sets $N^+(v_i)$, adjacency testing can be
performed in $O(d)$ time, since two vertices $v_i$ and $v_j$ are
adjacent if and only if $v_j\in N^+(v_i)$ where $i<j$; see
\citep{CE91}.


\begin{center}
\framebox{
  \begin{minipage}{120mm}
    \textsc{DegenerateCliques}($G$, $d$)\\
    \emph{input}: a $d$-degenerate  graph $G$\\
    \emph{output}: all cliques in $G$ \medskip \hrule \medskip
    \begin{tabular}{ll}
      1. &         compute a $d$-degenerate ordering $(v_1,\dots,v_n)$ of $G$\\
      2. &    compute the sets $\{N^+(v_i):1\leq i\leq n\}$\\
      3. &    for $i:=1,\dots,n$ do\\
      4. &    \hspace*{10mm} $\textsc{AllCliques}(G[\{v_i\}\cup N^+(v_i)]$\\
     5. &    end-for\\
    \end{tabular}
  \end{minipage}} 
\end{center}
\medskip

\begin{theorem}
  \thmlabel{DegenCliques} If $G$ is a $d$-degenerate $n$-vertex graph,
  then $\textsc{DegenerateCliques}(G,d)$ outputs all the cliques in
  $G$ in time $O(d\,2^dn)$.
\end{theorem}

\begin{proof} 
  If $C$ is a clique of $G[N^+(v_i)]$ then $C\cup\{v_i\}$ is a clique
  of $G$. Thus every set output by the algorithm is a clique of
  $G$. Conversely, if $S$ is a clique of $G$, and $i$ is the minimum
  integer such that $v_i\in S$, then $S\setminus\{v_i\}$ is a clique
  of $G[N^+(v_i)]$, and $S$ is output by the algorithm. Now consider
  the time complexity. Since adjacency testing can be performed in
  $O(d)$ time, the subgraph $G[\{v_i\}\cup N^+(v_i)]$ can be
  constructed in $O(d^3)$ time.  By \thmref{NonRecCliques}, the call
  to \textsc{AllCliques} takes $O(d\,2^d)$
  time. 
  Hence the total time is $O(d\,2^dn)$.
\end{proof}


Since $H$-minor free graphs are $O(t\sqrt{\log t})$-degenerate, where
$t=|V(H)|$, \thmref{DegenCliques} implies:

\begin{corollary}
  \corlabel{MinorCliques} For every fixed graph $H$, there is a linear
  time algorithm to list all cliques in a given $H$-minor-free
  graph.\qed
\end{corollary}

\section{Graph Minor Decomposition}


This section first describes the Robertson-Seymour decomposition
theorem characterising the structure of $H$-minor-free graphs, and
then describes the analogous decomposition theorem for odd-minor-free
graphs.  We need a number of definitions.

An \emph{embedding} refers to a 2-cell embedding of a graph in a
(orientable or non-orientable) surface; that is, a drawing of the
vertices and edges of the graph as points and arcs in the surface such
that every face (region outlined by edges) is homeomorphic to a disk;
see \cite{MoharThom}.

Let $I$ be a linearly ordered set.  A \emph{path decomposition} of a
graph $G$ is a sequence $(\bag_i:i\in I)$ of subsets of $V(G)$ called
\emph{bags} such that:
\begin{enumerate}
\item $\bigcup_{i\in I}\bag_i= V(G)$;
\item for each edge $uv \in E(G)$, there exists $i\in I$ such that
  both $u$ and $v$ are in $\bag_i$; and
\item for each vertex $v\in V(G)$, the set $\{i : v \in \bag_i\}$ is a
  sub-interval of $I$.
\end{enumerate}
The \emph{width} of $(\bag_i:i\in I)$ is the maximum cardinality of a
bag minus~$1$. The \emph{pathwidth} of a graph $G$ is the minimum
width over all possible path decompositions of $G$.


At a high level, the Robertson-Seymour decomposition theorem says that
for every graph $H$, every $H$-minor-free graph can be expressed as a
tree structure of pieces, where each piece is a graph that can be
drawn in a surface in which $H$ cannot be drawn, except for a bounded
number of ``apex'' vertices and a bounded number of local areas of
non-planarity called ``vortices''.  Here the bounds depend only
on~$H$. Each piece in the decomposition is ``$h$-almost-embeddable''
where $h$ is a constant depending on the excluded minor~$H$. Roughly
speaking, a graph $G$ is \emph{$h$-almost embeddable} in a surface
$\Sigma$ if there exists a set $A\subseteq V(G)$ of size at most $h$,
such that $G-A$ can be obtained from a graph embedded in $\Sigma$ by
attaching at most $h$ graphs of pathwidth at most $h$ to within $h$
faces in an orderly way. The elements of $A$ are called \emph{apex
  vertices}.

More precisely, , a graph $G$ is \emph{$h$-almost embeddable} in a
surface $\Sigma$ if there exists a set $A\subseteq V(G)$ of size at
most $h$ such that $G-A$ can be written $G_0\cup G_1\cup \cdots \cup
G_h$, where
\begin{itemize}
\item $G_0$ has an embedding in $\Sigma$;
\item the graphs $G_1,\dots,G_h$, called \emph{vortices}, are pairwise
  disjoint;
\item there are faces $F_1,\dots,F_h$ of $G_0$ in $\Sigma$, and there
  are pairwise disjoint disks $D_1,\dots, D_h$ in $\Sigma$, such that
  for each $i\in\{1,\dots,h\}$,
  \begin{itemize}
  \item $D_i\subset F_i$ and $U_i:=V(G_0)\cap V(G_i)= V(G_0)\cap D_i$;
    and
  \item if $U_i$ is linearly ordered around the boundary of $F_i$,
    then $G_i$ has a path decomposition $(\bag_u:u\in U_i)$ of width
    less than $h$, such that $u\in \bag_u$ for each $u\in U_i$.
  \end{itemize}
\end{itemize}

The pieces of the decomposition are combined according to
``clique-sum'' operations, a notion which goes back to the
characterisations of $K_{3,3}$-minor-free and $K_5$-minor-free graphs
by \citet{Wagner37}. Suppose $G_1$ and $G_2$ are graphs with disjoint
vertex sets and let $k \geq 0$ be an integer. For $i=1,2$, suppose
that $W_i\subseteq V(G_i)$ is a $k$-clique in $G_i$. Let $G_{i}'$ be
obtained from $G_{i}$ by deleting some (possibly no) edges from the
induced subgraph $G_{i}[W_i]$ with both endpoints in $W_{i}$.
Consider a bijection $h : W_1\rightarrow W_2$. A \emph{$k$-sum} $G$ of
$G_1$ and $G_2$, denoted by $G=G_{1}\oplus_{k} G_{2}$ or simply by
$G=G_{1}\oplus G_{2}$ is the graph obtained from the union of $G'_1$
and $G'_2$ by identifying $w$ with $h(w)$ for all $w\in W_1$.  A
\emph{$(\leq k)$}-sum is a $k'$-sum for some $k'\leq k$. Note that
$\oplus$ is not uniquely defined.

Now we can finally state a precise form of the decomposition theorem:
\begin{theorem} {\rm \cite[Theorem~1.3]{RS-GraphMinorsXVI-JCTB03}}
  \thmlabel{general} For every graph $H$, there exists an integer
  $h\geq 0$ depending only on $|V(H)|$ such that every $H$-minor-free
  graph can be obtained by $(\leq h)$-sums of graphs that are
  $h$-almost-embeddable in some surfaces in which $H$ cannot be
  embedded.
\end{theorem}
In particular, if $H$ is fixed then a surface in which $H$ cannot be
embedded has bounded Euler genus.  Thus the summands in
\thmref{general} are $h$-almost embeddable in surfaces of bounded
Euler genus. A graph is \emph{$h$-almost embeddable} if it is
$h$-almost embeddable in a surface of Euler genus at most $h$.

We now describe a decomposition theorem for odd-minor-free graphs by
\citet{DHK-SODA10}. This result generalises \thmref{general}. A graph
$G$ is \emph{$h$-almost bipartite} if $G-A$ is bipartite for some set
$A\subseteq V(G)$ with $|A|\leq h$.

\begin{theorem}[\citep{DHK-SODA10}]
  \thmlabel{OddStructure} For every fixed integer $t$, there is a
  constant $h$ such that every odd-$K_t$-minor-free graph~$G$ can be
  obtained by $(\leq h)$-sums of $h$-almost bipartite graphs and
  $h$-almost embeddable graphs.
\end{theorem}


\section{Listing Cliques in Odd-Minor-Free Graphs}

This section describes an algorithm for finding all the cliques in a
graph $G$ excluding a fixed odd-minor. The time complexity is
$O(n^3)$. Thus, we may assume that $G$ is represented by an adjacency
matrix (which takes $O(n^2)$ time to pre-compute), and adjacency
testing can be performed in $O(1)$ time.

\begin{lemma}
  \lemlabel{AlmostBipartite} Let $G$ be an $h$-almost-bipartite graph
  on $n$ vertices. Then $G$ contains at most $2^hn^2+2$ cliques.
\end{lemma}

\begin{proof}
  $G-A$ is bipartite for some $A\subseteq V(G)$ with $|A|\leq h$.
  Since $G-A$ is triangle-free, the cliques in $G-A$ are precisely
  $E(G-A)\cup V(G-A)\cup\{\emptyset\}$. There are at most
  $\frac{1}{4}(n-|A|)^2+n-|A|+1$ such cliques.
There are at most $2^{|A|}$ cliques in $G[A]$.
Every clique in $G$ is the union of a
  clique in $G-A$ and a clique in $G[A]$. Thus $G$ contains at most
  $2^{|A|}(\frac{1}{4} (n-|A|)^2+n-|A|+1)\leq 2^hn^2+2$ cliques.  
\end{proof}

\begin{lemma}
  \lemlabel{AlmostEmbeddable} Let $G$ be an $h$-almost embeddable
  graph on $n$ vertices.  Then, for some $h'$ and $h''$ that only
  depend on $h$, $G$ contains at most $h'n$ cliques, and they can be
  listed in $O(h''n)$ time.
\end{lemma}

\begin{proof}
  It is well known that $G$ contains no $K_{h'}$-minor, for some $h'$
  depending only on $h$ (see \citep{CliqueMinors} for a tight bound on
  $h'$). Thus $G$ is $O(h'\sqrt{\log h'})$-degenerate, and the claim
  follows from \corref{MinorCliques}.
\end{proof}

\begin{lemma}
  \lemlabel{BasicCliqueSum} Let $c>0$. Let $G$ be a $k$-sum of graphs
  $G_1$ and $G_2$, where each $G_i$ has $n_i$ vertices and contains at
  most $cn_i^2$ cliques. Assume that $n_1\geq\frac{k^2}{2}+k$ and $G$
  has $n$ vertices.  Then $G$ contains at most $cn_1^2+cn_2^2$
  cliques, which is at most $cn^2$.
\end{lemma}

\begin{proof}
  Since $n_1\geq \frac{k^2}{2}+k$ and $n_2\geq k+1$,
\begin{align*}
n_1\geq \frac{k^2}{2}+k\geq  \frac{k^2}{2(n_2-k)} + k
&=\frac{2k(n_2-k) + k^2}{2(n_2-k)} \\
&=\frac{k(2n_2- k)}{2(n_2-k)}
\enspace.
\end{align*}
Hence $2n_1n_2-2kn_1 \geq 2kn_2- k^2$, implying
\begin{align*}
n^2&=(n_1+n_2-k)^2\\
&=n_1^2+n_2^2+2n_1n_2-2kn_1-2kn_2+k^2 \\
&\geq n_1^2+n_2^2\enspace.
\end{align*}
Each clique in $G$ is a clique of $G_1$ or $G_2$. Thus $G$ contains at
most $cn_1^2+cn_2^2\leq cn^2$ cliques.
\end{proof}

\begin{lemma}
  \lemlabel{CliquesSum} Let $k$ be a positive integer.  Let
  $G_1,\dots,G_p$ be graphs, such that each $G_i$ has $n_i$ vertices
  and contains at most $f(k)\cdot n_i^2$ cliques, for some function
  $f$.  Furthermore, suppose that each $G_i$ contains no
  $k$-clique. Let $G$ be an $n$-vertex graph obtained by $(\leq
  k)$-sums of $G_1,\dots,G_p$. Then for some function $f'$ depending
  on $f$ and $k$, $G$ contains at most $f'(k)\cdot n^2$ cliques.
\end{lemma}

\begin{proof} 
  The construction of $G$ defines a binary tree $T$ rooted at some
  node $r$, and associated with each node $v$ of $T$ is a subgraph
  $G_v$ of $G$, such that $G_r=G$; $G_1,\dots,G_p$ are the subgraphs
  associated with the leaves of $T$; and $G_v=G_u\oplus_{\leq k} G_w$
  for each non-leaf node $v$ with children $u$ and $w$. Let $n_v$ be
  the number of vertices in each $G_v$.  Say $G_v$ is \emph{small} if
  $n_v<\frac{k^2}{2}+k$.

  Let $T'$ be the subtree of $T$ obtained by applying the following
  rule until it cannot be further applied: If $u$ and $w$ are leaf
  nodes with a common parent $v$, and both $G_u$ and $G_w$ are small,
  then delete $u$ and $w$. The remainder of the proof focuses on $T'$.

  We now prove (by induction, working from the leaves of $T'$ up
  through the tree) that each subgraph $G_v$ contains at most
  $f'(k)\cdot n_v^2$ cliques, where
  $f'(k):=\max\{f(k),2^{k^2+2k}\}$. If $v$ is a leaf of $T$ then this
  hypothesis holds by assumption. If $v$ is a leaf of $T'$ but not of
  $T$, then $G_u$ and $G_w$ are small, where $u$ and $w$ are the
  children of $v$ in $T$. In this case $n_v\leq k^2+2k$, implying
  $G_v$ contains at most $2^{k^2+2k}\leq f'(k)\cdot n_v^2$
  cliques. Thus the hypothesis again holds.

  Now consider a non-leaf node $v$ of $T'$. Let $u$ and $w$ be the
  children of $v$. We have $G_v=G_u\oplus_\ell G_w$ for some $\ell\leq
  k$.  By induction, $G_u$ contains at most $f'(k)\cdot n^2_u$
  cliques, and $G_w$ contains at most $f'(k)\cdot n^2_w$
  cliques. Suppose that $G_u$ and $G_w$ are both small. If $u$ and $w$
  are both leaves in $T$ then the above rule is applicable. Otherwise,
  without loss of generality, $w$ is not a leaf in $T$, in which case
  every descendent subgraph of $w$ is small, implying the subtree
  rooted at $w$ contains two leaves for which the above rule is
  applicable. Hence at least one of $G_u$ and $G_w$ is not small. Thus
  \lemref{BasicCliqueSum} is applicable with $c=f'(k)$.  Hence $G_v$
  contains at most $f'(k)\cdot n_u^2+f'(k)\cdot n_w^2$ cliques, which
  is at most $f'(k)\cdot n_v^2$ cliques.  In particular, $G=G_r$
  contains at most $f'(k)\cdot n^2$ cliques, as claimed. Observe that
  the above argument actually proves that the sum of $n_u^2$, taken
  over all leaf nodes $u$ in $T'$, is at most $n^2$.
%
%
\end{proof}

\begin{proof}[Proof of \thmref{OddCliques}]
  By \thmref{OddStructure}, $G$ is the $(\leq h)$-sum of graphs
  $G_1,\dots,G_p$, where each $G_i$ is $h$-almost bipartite or
  $h$-almost embeddable in a surface of Euler genus $h$.  By
  \twolemref{AlmostBipartite}{AlmostEmbeddable}, for some $h'$ that
  only depends on $h$, if each $G_i$ has $n_i$ vertices, then $G_i$
  contains at most $h'n_i^2$ cliques. Note that $G$ contains no
  $h$-clique.  By \lemref{CliquesSum}, $G$ contains at most $h''\,n^2$
  cliques, for some $h''$ depending only on $h$. 
By \thmref{NonRecCliques}, the cliques in
  $G$ can be output in $O(h''\,n^2)$ time by algorithm \textsc{AllCliques}($G$).
\end{proof}

Note that reference \citep{DHK-SODA10} describes a polynomial time
algorithm for computing the decomposition described in
\thmref{OddStructure}. However, by using \thmref{NonRecCliques} it suffices
to merely prove an upper bound on the number of cliques in an
odd-minor-free graph, to obtain an efficient algorithm for listing all
cliques. 
 

\subsection*{Acknowledgements} Many thanks to the referees of the
conference version of this paper who found some minor errors that have now been corrected. 


\def\soft#1{\leavevmode\setbox0=\hbox{h}\dimen7=\ht0\advance \dimen7
  by-1ex\relax\if t#1\relax\rlap{\raise.6\dimen7
  \hbox{\kern.3ex\char'47}}#1\relax\else\if T#1\relax
  \rlap{\raise.5\dimen7\hbox{\kern1.3ex\char'47}}#1\relax \else\if
  d#1\relax\rlap{\raise.5\dimen7\hbox{\kern.9ex \char'47}}#1\relax\else\if
  D#1\relax\rlap{\raise.5\dimen7 \hbox{\kern1.4ex\char'47}}#1\relax\else\if
  l#1\relax \rlap{\raise.5\dimen7\hbox{\kern.4ex\char'47}}#1\relax \else\if
  L#1\relax\rlap{\raise.5\dimen7\hbox{\kern.7ex
  \char'47}}#1\relax\else\message{accent \string\soft \space #1 not
  defined!}#1\relax\fi\fi\fi\fi\fi\fi}

\end{document}